\documentclass[reqno]{article}

\usepackage{xcolor}
\usepackage{amssymb}
\usepackage{amsfonts}
\usepackage{amsmath}
\numberwithin{equation}{section}
\usepackage{enumerate}
\newtheorem{theorem}{Theorem}[section]
\newtheorem{proposition}[theorem]{Proposition}
\newtheorem{lemma}[theorem]{Lemma}
\newtheorem{corollary}[theorem]{Corollary}
\newtheorem{definition}[theorem]{Definition}
\newtheorem{example}[theorem]{Example}

\title {Root functions of a meromorphic matrix function and applications} 
\author{Author: Muhamed Borogovac}

\begin{document}
\maketitle 
This paper is dedicated to the memory of Prof. Heinz Langer.
\\

\textbf{Abstract.} A practical method is presented for determining root and pole cancellation functions of a matrix function $Q(z)$ meromorphic on the extended complex plane $\bar{\mathbb{C}} := \mathbb{C} \cup { \infty }$. The method is applied to nonlinear systems of $n \in \mathbb{N}$ differential equations of order $m \in \mathbb{N}$ with $n$ unknown functions $u_{i}(t)$, where $i = 1, \ldots, n$.
 
For a function $Q \in \mathcal{N}_{\kappa}(\mathcal{H})$, $\kappa \in \mathbb{N} \cup \lbrace 0 \rbrace $, possessing a pole of order $m \in \mathbb{N}$ at infinity, the following representation is established:
\[
Q(z)=(z-\beta)^{m}\tilde{Q}(z), \, z\in \mathcal{D}(Q),
\]
where $\beta \in \mathbb{R}$ is a regular point of $Q$, and $\tilde{Q} \in \mathcal{N}_{\tilde{\kappa}}(\mathcal{H})$ is holomorphic at $\infty$. Unlike the Krein–Langer representation of $Q$, which involves a linear relation $A$, this representation employs a bounded operator $\tilde{A}$ in the Krein-Langer representation of $\tilde{Q}$. The operator $\tilde{A}$ and the relation $A$ have identical spectra, except at $\beta$ and $\infty$. We demonstrate how to obtain this representation for a given meromorphic function $Q \in N^{n \times n}_{\kappa}$ using the root functions developed in this work. 
\\
\\
MSC (2020) 47A56 34M04 47B50 46C20 
\\
\\
Keywords: Meromorphic matrix function; Root functions; Operator representation of a generalized Nevanlinna function; Nonlinear system of differential equations.

\section{Introduction}\label{s2}
\textbf{1.1.} Let $\mathbb{N}$, $\mathbb{R}$, and $\mathbb{C}$ denote the sets of positive integers, real numbers, 
and complex numbers, respectively. Let $(.,.)$ denote the (definite) scalar 
product in a Hilbert space $\mathcal{H}$, and let $\mathcal{L}(\mathcal{H})$ denote the 
Banach space of bounded linear operators on $\mathcal{H}$. We will mostly focus on the case $\mathcal{H} = \mathbb{C}^n$, where $n \in \mathbb{N}$, due to certain applications we have developed. Nevertheless, some statements will be proved in a more general setting.

We assume that $\mathcal{D}(Q) \subseteq \mathbb{C}$ is a domain, and that $Q:\mathcal{D}(Q) \to \mathcal{L}(\mathbb{C}^n)$ is a matrix-valued complex function with $\det Q(z) \not\equiv 0$.

To investigate the zeros and poles of complex matrix functions, it is necessary to define their so-called root and pole cancellation functions. The following definition is consistent with those given in \cite{GS,GLR,B1,B5}.   
 
\begin{definition}\label{definition42} Let $Q(z)$ be an $n\times n$ meromorphic matrix function on $\bar{\mathbb{C}}$. Suppose that the $n$-dimensional vector function $\varphi(z)$ satisfies the following conditions at $\alpha \in \mathbb{C}$:
\begin{enumerate}[(a)]
\item $\varphi(\alpha) \ne 0$;
\item $(Q(z)\varphi(z))^{(l)} \to 0$ as $z \to \alpha$, for $0 \le l \le m-1$, where $l, m-1 \in \mathbb{N} \cup \{0\}$.
\end{enumerate}
Then $\varphi(z)$ is called a \textit{root function} of $Q(z)$ of order at least $m$ at the critical point (or zero)~$\alpha$. If $m$ is the maximal number for which condition~(b) holds, then $\varphi(z)$ is said to be a root function of \textit{exact order} (or \textit{exact partial multiplicity}) $m$ at~$\alpha$.  

We call $\varphi(\alpha)$ an \textit{eigenvector corresponding to the eigenvalue}~$\alpha$. If $\varphi(z)$ is a root function of maximal order $k$ among all root functions with the same eigenvector $\varphi(\alpha)$, then $\varphi(z)$ is called the \textit{canonical root function} associated with $\varphi(\alpha)$, and the sum of all orders of the canonical root functions at~$\alpha$ is called the (total) multiplicity of the zero at~$\alpha$. The maximal order among all canonical root functions at~$\alpha$ is called the \textit{order of the zero (or critical point)} at~$\alpha$.

If the matrix function $Q(z)$ is invertible, i.e. $\chi(z) := \det Q(z) \not\equiv 0$, then $\varphi(z)$ is also called a \textit{pole function} of $Q^{-1}(z)$ of the same order. 
\end{definition}
The following definition is consistent with the definitions of pole functions and pole cancellation functions found in \cite{B1, B5, BLu, DLS}.

\begin{definition}\label{definition44} Let $Q(z)$ be an $n\times n$ meromorphic matrix function on $\bar{\mathbb{C}}$. Suppose that the n-dimensional vector functions $\hat{\varphi}(z)$  and $\psi(z)$ satisfy the following conditions at $\beta \in \mathbb{C}$:
\begin{enumerate}[(a)]
\item $\hat{\varphi} \left( z \right):=Q\left( z \right)\psi \left( z \right)\to c=:\hat{\varphi}(\beta) \ne 0,\, c\ne \infty ,\, as\, z\to \beta $, 
\item $\psi^{\left( l \right)}\left( z \right)\to 0,\, as\, z\to \beta ,\, 0\le l\le m-1$, where $l,m-1\in \mathbb{N}\cup \lbrace 0\rbrace$.
\end{enumerate}
Then $\beta \in \mathbb{C}$ is called a pole of $Q\left( z \right)$. The functions 
$\hat{\varphi} \left( z \right)$ and $\psi (z)$ are called a pole and pole cancellation 
functions, respectively, of $Q\left( z \right)$ at $\beta $ of order at least $m$. If $m$ is the 
maximal number with property (b) for a particular vector $\hat{\varphi} 
\left( \beta \right)$, then $\hat{\varphi} \left( z \right)$ and $\psi \left( z 
\right)$ are called the \textit{pole} and \textit{pole cancellation functions}, respectively, of $Q\left( z \right)$ at 
$\beta $ of \textit{exact order} $m$.

If $\psi (z)$ is a pole cancellation function of maximal order $k$ among all pole cancellation functions associated with the same limit $\hat{\varphi} (\beta)$, then $\psi (z)$ is referred to as the canonical pole cancellation function associated with the vector $\hat{\varphi} (\beta)$, and the sum of all orders of those canonical functions at $\beta$ is called the (total) multiplicity of the pole at $\beta$. The maximal order among all canonical pole cancellation functions at $\beta$ is called the order of the pole at $\beta$.

If the matrix function $Q\left( z \right)$ is invertible, then $\hat{\varphi}\left( z \right)$ is called a root function and $\psi \left( z \right)$ is called a zero function of $Q^{-1}\left( z \right)$. 
\end{definition}
By definition, $\infty$ is a pole (respectively, a zero) of $Q$ if and only if $0$ is a pole (respectively, a zero) of $g(z) = Q(1/z)$. Hence, all of the preceding definitions extend naturally to the case $\alpha = \infty$.

The relationship between a root function $\varphi (z)$ of $Q\left( z \right)$ at a critical point $\alpha \in \mathbb{C}$ and a pole cancellation function 
$\hat{\psi} (z)$ of $Q^{-1}\left( z \right)$ at the pole $\alpha $ is:
\[
Q\left( z \right)\varphi \left( z \right)=\hat{\psi} \left( z \right), \, \,\, z\in \mathcal{D}\left( Q \right).
\]
Hence, a point $\beta \in \bar{\mathbb{C}}$ is a zero (respectively, a pole) of $Q^{-1}$ if and only if it is a pole (respectively, a zero) of $Q$ of the same order and multiplicities.

Recall that an operator-valued complex function $Q:\mathcal{D}\left( Q \right)\to \mathcal{L}(\mathcal{H})$ 
belongs to the class of \textit{generalized Nevanlinna functions} $\mathcal{N}_{\kappa }(\mathcal{H})$ if it is meromorphic in 
$\mathbb{C} \backslash \mathbb{R}$, symmetric, i.e. ${Q\left( \bar{z} \right)}^{\ast 
}=Q\left( z \right),\, \, \forall z\in \mathcal{D}\left( Q \right)$, and the 
Nevanlinna kernel 
\[
\mathcal{N}_{Q}\left( z,w \right):=\frac{Q\left( z \right)-{Q\left( w \right)}^{\ast 
}}{z-\bar{w}};\, z,\, w\in \mathcal{D}(Q)\cap \mathbb{C}^{+},
\]
has $\kappa $ negative squares. 
\\
\\
\textbf{1.2.} Since a generalized Nevanlinna function $Q\left( z \right)$ admits a 
Krein-Langer operator representation (see \cite{KL3,KL2,HSW,BeLu}), such representations can be used to express root and pole cancellation functions in terms of the Jordan chains of the associated operator or relation $A$; 
see, e.g., \cite{B1,BLu}. These root and pole cancellation functions play a key role in characterizing 
the orders and multiplicities of the zeros and poles of a generalized 
Nevanlinna function $Q\left( z \right)$, as well as in identifying 
generalized poles of non-positive type; see, e.g., \cite{B1,BL,DLS,BLu}. 

Nonetheless, \textbf{the main difficulty remains in the explicit determination of these root and pole cancellation functions}. In particular, for a given generalized Nevanlinna function, it is often infeasible to explicitly determine the representing operator or relation $A$ in the corresponding Krein-Langer representation of $Q$. Without explicit knowledge of this operator or relation $A$ and the corresponding Jordan chains, one cannot determine the associated pole cancellation functions by means of the results presented in the cited papers.
\\

Currently, feasible methods exist only for finding root functions of matrix polynomials; see \cite{B5}. \textbf{One of the main objectives of this paper is to extend the method introduced in \cite{B5} by constructing explicit root and pole cancellation functions for a broader class of functions—specifically, meromorphic matrix functions on the extended complex plane.} Notably, our analysis is not restricted to the generalized Nevanlinna class or to Krein–Langer representations.

The new method for determining root and pole cancellation functions is presented in Section \ref{s4}, in Theorems \ref{theorem48} and \ref{theorem412}, Proposition \ref{proposition414}, and Example \ref{example418}.

In Section \ref{s6}, we demonstrate the advantages of the method developed in Section \ref{s4} over the logarithmic residue theorem for determining the orders and multiplicities of zeros and poles of meromorphic matrix functions on $\bar{\mathbb{C}}$.

In Section \ref{s8}, we study the general form of the nonlinear system of differential equations
\begin{equation}
\label{eq22} 
\begin{array}{*{20}c}
\frac{a_{11}^{1}}{u'_{1}}+\frac{a_{12}^{1}}{u'_{2}}+\frac{a_{11}^{0}}{u_{1}}+\frac{a_{12}^{0}}{u_{2}}=0,\\
\frac{a_{21}^{1}}{u'_{1}}+\frac{a_{22}^{1}}{u'_{2}}+\frac{a_{21}^{0}}{u_{1}}+\frac{a_{22}^{0}}{u_{2}}=0,\\
\end{array} 
\end{equation}
where $a_{ij}^{k}$ ($ k=0,1; \, i,j=1,2$) are constants, and  $u_{1}(t)$ and $u_{2}(t)$ are unknown functions. 

Currently, no general method exists for solving such nonlinear systems. We consider a generalization of system~\eqref{eq22} to a system of $n \in \mathbb{N}$ differential equations of order $m \in \mathbb{N}$ with $n$ unknown functions $u_{j}$, $j = 1, \ldots, n$. By constructing a meromorphic matrix function associated with the system and applying the results from Section~\ref{s4}, we obtain its solutions. The solutions are given in Theorem \ref{theorem82}, and the method is illustrated in Example \ref{example84}. 

Assume now that $\infty \in \bar{\mathbb{C}} $ is a pole of a function $Q\in \mathcal{N}_{\kappa}(\mathcal{H})$. Then $Q\left( z \right)$ must be represented by a self-adjoint 
linear relation $A$ in the Krein–Langer type representation; see, for example, \cite[p.114]{DLS} or \cite[Theorem 4.2]{HSW}. 

In Theorem \ref{theorem104}, we establish  the following factorization of $Q\in \mathcal{N}_{\kappa}(\mathcal{H})$:
\begin{equation}
\label{eq24} 
Q(z)=(z-\beta)^{m}\tilde{Q}(z), \, z\in \mathcal{D}(Q),
\end{equation}
where $\beta \in \mathbb{R}$ is a regular point of $Q$, $m\in \mathbb{N}$ is the order of the pole at $\infty$, and the function $\tilde{Q}\in \mathcal{N}_{\tilde{\kappa}}(\mathcal{H}),\, \tilde{\kappa}\in \lbrace 0 \rbrace \cup \mathbb{N}$, is holomorphic at $\infty $. Via the Krein–Langer representation of $\tilde{Q}$, \textbf{the function $Q$ can now be represented by a well-behaved bounded self-adjoint operator $\tilde{A}$, rather than by a linear relation $A$. The operator $\tilde{A}$ has the same spectrum as the relation $A$, except that the eigenvalue at infinity of $A$ is replaced by a finite real eigenvalue $\beta$ of $\tilde{A}$ of the same order.} In Example \ref{example105}, it is shown that, in general, the negative indices and multiplicities at $\beta$ and $\infty$ need not be the same, whereas Theorem \ref{theorem104} provides the relation between the negative indices.

Moreover, if $Q$ is a meromorphic matrix function with real poles, then the method developed in Section \ref{s4} can be used to easily determine the poles and pole cancellation functions of $Q$. \textbf{This enables the explicit construction of the bounded representing operator $\tilde{A}$ and the realization (\ref{eq24}) of $Q$.} This construction is demonstrated in Example \ref{example106}.

We note that a different factorization of operator-valued generalized Nevanlinna functions was studied in \cite{Lu2}.

\section{Finding canonical sets of root and pole cancellation functions}\label{s4}
\textbf{2.1.} The following lemma summarizes considerations in \cite[Section 2.1]{B5}.

\begin{lemma}\label{lemma46} Let $L(z)$ be an $n\times n$ complex matrix polynomial 
with characteristic polynomial $\chi \left( z \right):=\det {L\left( z 
\right)}\not\equiv 0$. Then there exist an $n\times n$ diagonal matrix 
polynomial $D\left( z \right)$ and two $n\times n$ invertible 
matrix polynomials $S\left( z \right)$ and $T\left( z \right)$, each with 
constant (nonzero) determinants, such that
\begin{equation}
\label{eq42}
D\left( z \right)\mathrm{=}S\left( z \right)L\left( z \right)T\left( z 
\right).
\end{equation}
Moreover, the polynomial
\[
D\left( z \right)=diag\left[ d_{1}\left( z \right),... , d_{n}\left( z \right) \right]
:=\left( {\begin{array}{*{20}c}
d_{1}\left( z \right) & \mathellipsis  & 0\\
\mathellipsis  & \mathellipsis  & \mathellipsis \\
0 &\mathellipsis  & d_{n}\left( z \right)\\
\end{array} } \right),
\]
is uniquely determined up to the lineup and 
constant multiples of the scalar polynomials $d_{\mathrm{i}}\left( z \right)$ on 
the main diagonal, depending on the choice of the 
matrices of elementary transformations $S(z)$ and $T(z)$.
\end{lemma}

\begin{definition}\label{definition47}  A function $Q \in \mathcal{L}(\mathcal{H})$, where $\mathcal{H}$ is a Hilbert space, is said to be convergent at 
$\infty $ with limit $S$ if 
\[
\lim\limits_{z\to \infty}{{Q\left( z \right)}}=S
\]
In this case we define $Q(\infty):=S$. 

The function $Q$ is sad to be holographic at $\infty$ if the following limit exists:
\[
Q'\left( \infty \right):=\lim\limits_{z\to \infty}{z\left( Q\left( z \right)-S \right)},
\]
\end{definition}
for some $S\in \mathbb{C}$. The following theorem generalizes Lemma \ref{lemma46}.
\begin{theorem}\label{theorem48} Let $Q(z)$ be an $n\times n$ meromorphic matrix function 
on $\bar{\mathbb{C}}$ with characteristic function $\chi \left( z \right):=\det 
{Q\left( z \right)}\not\equiv 0$. Then there exist a diagonal $n\times n$ 
matrix function $\tilde{D}\left( z \right)$ and two invertible $n\mathrm{\times }n$ 
matrix polynomials $S\left( z \right)$ and $T\left( z \right)$, 
each with constant determinants, such that
\begin{equation}
\label{eq48}
\tilde{D}\left( z \right)=S\left( z \right)Q\left( z \right)T\left( 
z \right),
\end{equation}
where
\begin{equation}
\label{eq410}
\tilde{D}\left( z \right)=\operatorname{diag}\left[ 
\tilde{d}_{1}\left( z \right),...\, ,\tilde{d}_{n}\left( z \right) 
\right],
\end{equation}
\[
\tilde{d}_{i}\left( z \right)=\frac{\tilde{p}_{i}\left( z 
\right)}{\tilde{q}_{i}\left( z \right)}\mathrm{,\, }i=1,\, 
\mathellipsis ,\, n,
\]
with $\tilde{p}_{i}\left( z \right)$ and $\tilde{q}_{i}\left( z \right)$ being scalar polynomials. 

The function $\tilde{D}(z)$ is uniquely determined up to the lineup and 
constant multiples of the functions $\tilde{d}_{\mathrm{i}}\left( z \right)$ 
on the main diagonal.
\end{theorem}
\textbf{Proof}.
For a matrix-valued function $Q(z)$ that is meromorphic on $\bar{\mathbb{C}}$, the point $\infty$ is either a 
point of holomorphy or a pole, that is, an isolated singularity. Hence, there exists a disc centered at $z\mathrm{=0}$ with sufficiently large 
radius $R>0$ such that $Q\left( z \right)$ is holomorphic outside this disc, 
except possibly at $\infty $. Since the disc is bounded and all 
singularities of $Q$ are poles, the Bolzano–Weierstrass theorem implies that there are only finitely many poles, say 
$\beta_{1},\, \beta_{2},\, 
\mathellipsis ,\, \beta_{l}$ with respective orders 
$t_{1},\, t_{2}\mathrm{,\, \mathellipsis ,\, 
}t_{l}$, where $\beta_{i}\ne \beta_{j}$ 
for $i\ne j$. Then the function
\begin{equation}
\label{eq414}
L\left( z \right)=\left( z-\beta_{1} 
\right)^{t_{1}}\mathellipsis \left( z-\beta_{l} 
\right)^{t_{l}}Q\left( z \right),
\end{equation}
is holomorphic on $\mathbb{C}$ and has a pole at $\infty $. It is 
well known that the only functions holomorphic on $\mathbb{C}$ with a pole at 
infinity are polynomials; hence $L(z)$ is a polynomial. Let us 
denote 
\[
q\left( z \right)=\left( z-\beta_{1} \right)^{t_{i}}\mathellipsis \left( z-\beta_{l} \right)^{t_{l}}.
\]
By Lemma \ref{lemma46},
\[
D\left( z \right)\mathrm{=}S\left( z \right)q\left( z \right)Q\left( z \right)T\left( z \right),
\]
where the matrices of elementary transformations $S(z)$ and $T(z)$ have constants determinants.

Dividing both sides by $q\left( z \right)$ yields (\ref{eq48}), where
\begin{equation}
\label{eq416}
\tilde{D}\left( z \right)=\frac{1}{q\left( z \right)}diag\, \left[ 
d_{1}\left( z \right),...\, ,d_{n}\left( z \right) \right].
\end{equation}
Since, $d_{i}\left( z \right)$ and $q\left( z \right)$ may have common 
factors that cancel, in general $\tilde{p}_{i}\left( z \right)\ne 
d_{i}\left( z \right)$ and $\tilde{q}_{i}\left( z \right)\ne q\left( z 
\right)$. 

The statement on the uniqueness of $\tilde{D}(z)$ follows directly from the corresponding result in Lemma \ref{lemma46} and from equation (\ref{eq416}).\hfill $\square$  

Note that identity (\ref{eq42}) in Lemma \ref{lemma46} is based on the Smith algorithm presented in \cite[Theorem S1.1]{GLR}. That algorithm produces the matrix $D(z)$ with additional 
properties, known as the \textit{Smith form} of the polynomial $L(z)$. However, since the Smith form is not required in our context, the lengthy procedure can often be avoided. A more efficient sequence of elementary transformations is often sufficient to diagonalize the polynomial $L(z)$ given by (\ref{eq414}), see Examples \ref{example418} and \ref{example84}.

The following corollary is a direct consequence of relation (\ref{eq414}).

\begin{corollary}\label{corollary49}
A matrix-valued function  $Q(z)$ is meromorphic on $\bar{\mathbb{C}}$ if and only if it can be written as the ratio of a matrix polynomial to a scalar polynomial.
\end{corollary}

Let us now arbitrarily select some $\tilde{d}_{i}\left( z \right),\, \, 
i=1,\mathellipsis ,n$. Suppose we find all $\gamma_{i}$ zeros $\alpha 
_{i,j},\, j=1,\, \mathellipsis ,\, \gamma_{i}$ of the equation 
$\tilde{d}_{i}\left( z \right)=0$, along with their orders $s_{ij}$, and all 
$\delta_{i}$ poles $\beta_{i,j},\, j=1,\, \mathellipsis ,\, \delta_{i}$, 
along with their orders$\, t_{ij}$, then
\begin{equation}
\label{eq418}
\tilde{d}_{i}\left( z \right)\, =\frac{\left( z-\alpha_{i,1} 
\right)^{s_{i1}}\mathellipsis \left( z-\alpha_{i,\gamma_{i}} 
\right)^{s_{{i\gamma }_{i}}}}{\left( z-\beta_{i,1} \right)^{t_{i1}}\mathellipsis \left( z-\beta_{i,\delta_{i}} \right)^{t_{{i\delta 
}_{i}}}}.
\end{equation}
In the sequel, $\sigma_{0}\left( Q \right)$ denotes the \textit{set of all finite zeros }and $\sigma 
_{p}\left( Q \right)$ denotes \textit{the set of all finite poles of} $Q(z)$. 

If $\sigma_{0}\left( Q \right)\cap \sigma_{p}\left( Q \right)\ne \emptyset$, then the poles and zeros of $Q(z)$ and $\chi (z)$ may not coincide, since any common ones cancel each other out in $\chi (z)$.

Now, fix $\alpha \in \sigma_{0}\left( Q \right)$ and $\beta \in \sigma_{p}\left( Q \right)$.
According to (\ref{eq418}), a function $\tilde{d}_{i}\left( z \right)$ may have multiple zeros and multiple poles. Moreover, $\alpha \in \sigma_{0}(Q)$ may be a 
zero and $\beta \in \sigma_{p}\left( Q \right)$ may be a 
pole of several functions $\tilde{d}_{i}\left( z \right)$ in the matrix 
$\tilde{D}(z)$ defined by (\ref{eq410}). Therefore, for the matrix $\tilde{D}(z)$ we 
introduce the following index sets: 
\[
\Omega_{0}\left( \alpha \right):=\left\{ i\in \left\{ 1,\, \mathellipsis 
,n\, \right\}:\, \tilde{d}_{i}\left( \alpha \right)=0 \right\},
\]
\[
\Omega_{p}\left( \beta \right):=\left\{ i\in \left\{ 1,\, \mathellipsis ,n\, \right\}:\, {\tilde{d}_{i}\left( \beta \right)}^{-1}=0 \right\}.
\]
In the following theorem, the $i^{th}$ column vector of the matrix $T(z)$ is denoted by 
$T_{i}(z)$, and the $i^{th}$ column vector of the matrix $S^{-1}(z)$ is denoted by 
$\hat{S}_{i}(z)$. 

Unlike a matrix polynomial $L(z)$, for which $L(\alpha)$ and $\ker L(\alpha)$ are well defined, a meromorphic matrix function $Q(z)$ may be undefined at a critical point $\alpha$. In the following theorem, we use a generalized notion of the kernel at $\alpha$, denoted by $\ker Q\left( z \to \alpha \right)$, which remains meaningful when $\alpha \in \bar{\mathbb{C}}$ is both a zero and a pole of $Q$. The theorem below generalizes \cite[Theorem 2.4]{B5} from matrix polynomials to matrix functions meromorphic on $\bar{\mathbb{C}}$.

\begin{theorem}\label{theorem412} Let $Q(z)$ be an $n\times n$ meromorphic matrix 
function on $\bar{\mathbb{C}}$ with characteristic function $\chi \left( z 
\right):=\det {Q\left( z \right)} \not\equiv 0$. Assume that $\alpha \in 
\sigma_{0}\left( Q \right)$ is a zero of $\tilde{d}_{i}\left( z \right)$ of 
order $s_{i\alpha }$, and $\beta \in \sigma_{p}\left( Q \right)$ is a pole of 
$\tilde{d}_{j}\left( z \right)$ of order $t_{j\beta }$. Then: 

\begin{enumerate}[(i)]
\item The vector functions 
\begin{equation}
\label{eq424}
\varphi_{i}\left( z \right):=T_{i}\left( z \right),\, \, i\in \Omega_{0}\left( \alpha \right),
\end{equation}
are root functions of $Q\left( z \right)$ of orders $s_{i\alpha }$ corresponding to the eigenvectors $\varphi_{i}\left( \alpha \right)=T_{i}\left( \alpha \right)$; that is, they are the canonical 
root functions at $\alpha $.
\item For $i\in \Omega_{0}\left( \alpha \right)$, the vectors $\varphi_{i}\left( \alpha \right)$ are linearly independent, they determine the geometric multiplicity at $\alpha$ and the \textbf{generalized kernel of $Q$ at $\alpha$} by
\[
\ker {Q\left( z\rightarrow \alpha \right):=span\, \left\{ T_{i}\left( \alpha \right):\, 
i\in \Omega_{0}\left( \alpha \right) \right\}}.
\]
\item The functions 
\[
\psi_{j}\left( z \right):=T_{j}\left( z \right){\tilde{d}_{j}\left( z 
\right)}^{-1},\, \, j\in \Omega_{p}\left( \beta \right),
\]
are pole cancellation functions of the function $Q\left( z \right)$ of orders 
$t_{j\beta }$.
\item The vector functions 
\[
\hat{\varphi }_{j}\left( z \right):=\hat{S}_{j}\left( z \right),\, j\in 
\Omega_{p}\left( \beta \right),
\]
are root functions of $Q^{-1}\left( z \right)$ of orders $t_{j\beta }$, 
corresponding to the eigenvectors $\hat{\varphi }_{j}\left( 
\beta \right)=\hat{S}_{j}\left( \beta \right)$ ; that is, they are the 
canonical pole functions of $Q(z)$ at $\beta $.
\item For $j\in \Omega_{p}\left( \beta \right)$, the vectors $\hat{\varphi }_{j}\left( \beta \right)$, are linearly independent, they determine the geometric multiplicity at $\beta$ and the generalized kernel of $Q^{-1}$ at $\beta$ by
\[
\ker {Q^{-1}\left( z\rightarrow \beta \right)=span\, \left\{ \hat{S}_{j}\left( \beta 
\right):\, j\in \Omega_{p}\left( \beta \right) \right\}.}
\]
\item The functions $\hat{\psi }_{i}\left( z \right):=\hat{S}_{i}\left( z \right)\tilde{d}_{i}\left( z \right)$, for $i\in \Omega_{0}\left( \alpha \right)$, are \textbf{pole cancellation functions} of the inverse function $Q^{-1}\left( z \right)$ at the pole $\alpha $ of $Q^{-1}$, of orders $s_{i\alpha }$.
\end{enumerate}
\end{theorem}
\textbf{Proof.}

(i) Observe that $\det {S(z)}$ and $\det {T(z)}$ are constant nonzero 
values, since they are products of the determinants of elementary 
transformations. From identity (\ref{eq48}), it follows that
\begin{equation}
\label{eq428}
Q\left( z \right)T\left( z \right)={S\left( z \right)}^{-1}\tilde{D}\left( z 
\right), 
\end{equation}
Here, ${S\left( z \right)}^{-1}$ is also a matrix polynomial. Indeed, since 
$S\left( z \right)$ is a matrix polynomial, the elements of ${S\left( z 
\right)}^{-1}$ are given by the cofactors (i.e., minors) of $S\left( z 
\right)$ divided by the constant $det\, S\left( z \right)$; therefore, they 
are also polynomials.

Now, fix $\alpha \in \sigma_{0}\left( Q \right)$, and define the 
n-dimensional vector $P_{i}\left( \alpha \right)$ as follows:
\[
P_{i}\left( \alpha \right)=\left[ p_{1}\left( \alpha \right),\mathellipsis 
,\, p_{i}\left( \alpha \right),\, \mathellipsis ,\, p_{n}\left( \alpha 
\right) \right]^{T};\, p_{i}\left( \alpha \right)=\left\{ 
{\begin{array}{*{20}c}
1,\, if\, \tilde{d}_{i}\left( \alpha \right)=0, \\
0,\, othervise.\\
\end{array} } \right.
\]
By applying (\ref{eq428}) to $P_{i}\left( \alpha \right)$, we obtain: 
\[
Q\left( z \right)T\left( z \right)P_{i}\left( \alpha \right)={S\left( z 
\right)}^{-1}\tilde{D}\left( z \right)P_{i}\left( \alpha \right),\, \, \, \, 
\forall i\in \Omega_{0}\left( \alpha \right).
\]
Substituting $\tilde{D}(z)$ from (\ref{eq410}), we get:
\begin{equation}
\label{eq430}
Q\left( z \right)T_{i}\left( z \right)=\hat{S}_{i}\left( z 
\right)\tilde{d}_{i}\left( z \right),\, \, \, \forall i\in \Omega_{0}\left( 
\alpha \right).
\end{equation}
We are now ready to prove that the function $\varphi_{i}\left( z \right)$, 
as defined by (\ref{eq424}), is a root function of $Q$; that is, it satisfies 
conditions (a) and (b) of Definition \ref{definition44}.

(a) Since $det\, T(z)$ is a nonzero constant, the vectors $T_{i}\left( z \right)\ne 0$ for 
every $z\in \mathbb{C}$ and  $ i = 1, \ldots, n$, and they are linearly independent. For the selected eigenvalue 
$\alpha $, we have
\[
\varphi_{i}\left( \alpha \right)=T_{i}\left( \alpha \right)\ne 0,\, \, \, 
\forall i\in \Omega_{0}\left( \alpha \right).
\]
Thus, condition (a) is satisfied.

(b) Choose an arbitrary index $i\in \left\{ 1, \ldots, n
\right\}$ and let $z=\alpha $ be a zero of $\tilde{d}_{i}\left( z \right)$ 
of order $s_{i}$. Then we can write $\tilde{d}_{i}\left( z \right)\, =\, 
{(z\, -\, \alpha )}^{s_{i}}\, \hat{d}_{i}(z)$, where $\hat{d}_{i}(\alpha )\ne 0$. 
Using this, together with (\ref{eq430}), we obtain
\begin{equation}
\label{eq432}
Q\left( z \right)\varphi_{i}\left( z \right)={(z\, -\, \alpha )}^{s_{i}}\, 
\hat{S}_{i}\left( z \right)\hat{d}_{i}\left( z \right).
\end{equation}
Since both $\hat{S}_{i}\left( z \right)$ and $\hat{d}_{i}(z)$ are holomorphic at $z=\alpha $, it follows that
\[
\left( Q\left( z \right)\varphi_{i} \left( z \right) \right)^{\left( l 
\right)}\to 0,\, as\, z\to \alpha ,\, 0\le l \le s_{i}-1.
\]
Hence, $\varphi_{i}\left( z \right)$ satisfies condition $(b)$ of 
Definition \ref{definition44}, and is therefore a root function of $Q(z)$ at 
$\alpha $ of order at least $s_{i}$.

It remains to show that $s_{i}$ is the maximal order satisfying 
Definition \ref{definition44} (b). Assume to the contrary, that
\[
\left( Q\left( z \right)\varphi_{i} \left( z \right) \right)^{{(s}_{i})}\to 0,\, 
as\, z\to \alpha.
\]
Differentiating the right-hand side of equation (\ref{eq432}) $s_{i}$ times yields
\[
\left( {(z\, -\, \alpha )}^{s_{i}}\, \hat{S}_{i}\left( z 
\right)\hat{d}_{i}\left( z \right) \right)^{{(s}_{i})}=\left( z\, -\, \alpha 
\right)H\left( z \right)+\, \hat{S}_{i}\left( z \right)\hat{d}_{i}\left( z 
\right),
\]
where $H(z)$ is a function holomorphic at $\alpha$. Taking the 
limit as $z\to \alpha $, we obtain:
\[
\left( Q\left( z \right)\varphi_{i} \left( z \right) \right)^{{(s}_{i})}\to \, 
\hat{S}_{i}\left( \alpha \right)\hat{d}_{i}\left( \alpha \right)=0.\, 
\]
Since $\hat{d}_{i}\left( \alpha \right)\ne 0$, it must be $\hat{S}_{i}\left( 
\alpha \right)=0$. However, this contradicts the fact that $ 
det\, S\left( z \right)=\, c\ne \, 0$. This proves that 
Definition \ref{definition44} is satisfied with $s_{i}$ as the maximal order.

(ii) The vectors $\varphi_{i}\left( \alpha \right)$, for $i\in \Omega 
_{0}\left( \alpha \right)$, are linearly independent eigenvectors of $Q\left(z \right)$ at  $\alpha$, since they are the columns of 
the invertible matrix $T\left( \alpha \right)$. By the definition of 
$\Omega_{0}\left( \alpha \right)$, the set $\ker {Q\left( z\rightarrow \alpha \right)}$ contains all eigenvectors at $\alpha$, that is, $\dim  \ker {Q\left( z\rightarrow \alpha \right)}$ equals the geometric multiplicity at $\alpha$. Therefore, it is justified to refer to $\ker {Q\left( z\rightarrow \alpha \right)}$ as \textit{kernel}.

(iii), (iv) Clearly
\[
\tilde{D}^{-1}\left( z \right)=\left( {\begin{array}{*{20}c}
{\tilde{d}_{1}\left( z \right)}^{-1} & \mathellipsis & 
0\\
\mathellipsis & \mathellipsis  & \mathellipsis 
\\
0 & \mathellipsis & {\tilde{d}_{n}\left( z 
\right)}^{-1}\\
\end{array} } \right). 
\]
Then
\[
\left( \ref{eq48} \right)\Rightarrow S^{-1}\left( z 
\right)=Q\left( z \right)T\left( z \right){\tilde{D}\left( z \right)}^{-1}\Rightarrow 
\hat{S}_{j}\left( z \right)=Q\left( z \right)T_{j}\left( z 
\right){\tilde{d}_{j}\left( z \right)}^{-1}, \, j=1,...,n.
\]
This means that $\psi_{j}\left( z \right):=T_{j}\left( z 
\right){\tilde{d}_{j}\left( z \right)}^{-1}, \,j\in \Omega_{p}\left( \beta 
\right)$, is a pole cancellation function of $Q\left( z \right)$ of order $t_{j\beta}$, while 
$\hat{S}_{j}\left( z \right)$ is a root function of $Q^{-1}\left( z \right)$, 
i.e. it is a pole function of $Q(z)$ of order $t_{j\beta}$.

Statements (v) and (vi) follow directly from (ii) and (iii) by replacing 
$Q(z)$ by $Q^{-1}\left( z \right)$, $\tilde{D}\left( z \right)$ by 
$\tilde{D}^{-1}\left( z \right)$, and $T\left( z \right)$ by $S^{-1}\left( z 
\right)$. \hfill $\square$  
\\

Note that the functions $\varphi_{i}\left( z \right)$ and $\psi_{i}\left( z \right)$, for $i\in \Omega_{0}\left( \alpha \right)$, and the functions $\hat{\varphi }_{j}\left( z \right)$ and $\hat{\psi }_{j}\left( z \right)$, for $j\in \Omega_{p}\left( \beta\right)$, are not uniquely determined, as they depend on the choice of elementary transformations. However, according to Theorem \ref{theorem48}, they are canonical root and pole cancellation functions, and their orders are $s_{i\alpha }$ and $t_{jb}$, respectively. 

Therefore, we can now give the following equivalent definition of the orders and multiplicities of the zeros and poles of the meromorphic matrix function $Q$ on $\bar{\mathbb{C}}$, in terms of the diagonal matrix function $\tilde{D}(z)$. 

\begin{definition}\label{definition410}. If $\alpha \in \sigma_{0}\left( Q \right)$, then 
the number of the functions $\tilde{d}_{i}\left( z \right)$ that satisfy 
$\tilde{d}_{i}\left( \alpha \right)=0$ is called the geometric multiplicity 
of the zero $\alpha $. Let $s_{i}(\alpha )$ denote the order of the zero 
$\alpha $ of the function $\tilde{d}_{i}\left( \alpha \right)$. The numbers 
$s_{i}(\alpha )$ are called \textit{partial multiplicities} of the zero $\alpha $ of $Q(z)$. The largest among them is referred to as the order or rank of the zero $\alpha $, and the sum of all $s_{i}(\alpha )$ is called the \textbf{algebraic or total multiplicity} of the zero $\alpha $. 

If $\det {Q\left( z \right)}\not\equiv 0$, then for $\beta \in \sigma 
_{p}\left( Q \right)$, the geometric multiplicity, partial multiplicities, order 
and algebraic multiplicity of the pole $\beta $ of $Q(z)$ are defined as 
geometric multiplicity, partial multiplicities, order and algebraic multiplicity of 
$\beta $ as a zero of $Q^{-1}(z)$, respectively.
\end{definition}
\begin{proposition}\label{proposition414}
Let $\alpha \in \sigma_{0}\left( Q \right)$ be a zero (but not a pole) of the meromorphic matrix function $Q(z)$ on $\bar{\mathbb{C}}$.
A function $\varphi(z)$ is a root function of $Q(z)$ at $\alpha$ if and only if it is a root function of the same order of the matrix polynomial $L(z)$ given by~(\ref{eq414}). Therefore, the root functions of meromorphic matrix functions are essentially matrix polynomials.
\end{proposition}

\textbf{Proof.} Since $\alpha \in \sigma_{0}\left( Q \right)$ is not a pole, according to (\ref{eq414}) there exist a scalar polynomial $q(z)$ and a matrix polynomial $L(z)$ such that
\[
L(z)=q(z)Q(z)\, \wedge \, q(\alpha)\neq 0.
\]
Since $q(\alpha) \neq 0$, the function $\varphi(z)$ is a root function of $Q(z)$ if and only if it is a root function of $L(z)$ of the same order.\hfill $\square$

Evidently, a root function $\varphi(z)$ satisfies all the conditions stated in Definition~\ref{definition42} to admit a Taylor expansion
\begin{equation}
\label{eq440}
\varphi(z) = \sum_{j=0}^{l-1} (z - \alpha)^{j} \varphi_{j}(\alpha) + (z - \alpha)^{l} \tilde{\varphi}(z).
\end{equation}
of any order $l \le m - 1$, where $l, m \in \mathbb{N} \cup \lbrace 0 \rbrace$ are as defined in that definition. 


If $L(z) = A - zI$ is a matrix polynomial with an $n \times n$ matrix $A$, then the vectors $\varphi_{0}(\alpha), , \varphi_{1}(\alpha), , \ldots, , \varphi_{l-1}(\alpha)$ in the expansion~(\ref{eq440}) satisfy
\begin{equation}
\label{eq437}
\varphi_{j-1}(\alpha) = (A - \alpha I)\varphi_{j}(\alpha), \qquad j = 1, , \ldots, , l - 1,
\end{equation}
that is, they constitute a Jordan chain of the matrix $A$ corresponding to the eigenvalue $\alpha$.

If $k = l$ is the maximal number for which equations~(\ref{eq437}) hold, then the set of vectors $\varphi_{0}(\alpha), , \varphi_{1}(\alpha), , \ldots, , \varphi_{k-1}(\alpha)$ is called a \textit{maximal} or \textit{canonical Jordan chain}, and the corresponding function $\varphi(z)$ is called a \textit{canonical root function}.

Analogously, if $Q(z)$ is a meromorphic matrix function, the vectors $\varphi_{0}(\alpha), , \ldots, , \varphi_{l-1}(\alpha)$ in~(\ref{eq440}) can be called \textit{generalized Jordan vectors}.

If $L(z) = A - zI$ and the root function $\varphi(z)$ is constructed using $l < k$ terms in the sum~(\ref{eq440}), then $\varphi(z)$ has exact order $l < k$. However, there are examples that this is not necessarily the case for the Taylor expansions~(\ref{eq440}) of the root functions of other polynomial or meromorphic matrix functions.
\\
\\
\textbf{2.2.} All above claims are demonstrated in the following example. 

\begin{example}\label{example418} For the function
\[
Q\left( z \right)=\left( {\begin{array}{*{20}c}
\frac{z^{2}+4z+5}{z^{2}} &      -\frac{5}{z^{2}}\\
             &            \\
-\frac{5}{z^{2}} & \frac{25}{4z^{2}}\left( z^{2}+1 \right)\\
\end{array} } \right)\in \mathcal{N}_{\kappa}^{2\times2},
\]
use the method developed in Theorems~\ref{theorem48} and~\ref{theorem412} to find:
\begin{enumerate}[(i)]
\item eigenvalues and poles of $Q$;
\item the root functions of $Q$;
\item the pole cancellation functions of $Q$.
\end{enumerate}
\end{example}

\textbf{Solution.}

(i) and (ii) 

According to Theorem \ref{theorem48}, we must first obtain the representation (\ref{eq48}).
\[
q\left( z \right):=z^{2}\Rightarrow L\left( z \right)\mathrm{:=q(z)}Q\left( z 
\right)=\left( {\begin{array}{*{20}c}
z^{2}+4z+5 & -5\\
-5 & \frac{25}{4}\left( z^{2}+1 \right)\\
\end{array} } \right).
\]
Then
\[
\left( {\begin{array}{*{20}c}
I_{11} & L(z)\\
0 & I_{22}\\
\end{array} } \right)=\left( {\begin{array}{*{20}c}
1 & 0 &  z^{2}+4z+5 & -5\\
0 & 1 &  -5 & \frac{25}{4}\left( z^{2}+1\right) \\
0 & 0 &  1 & 0\\
0 & 0 &  0 & 1\\
\end{array} } \right)
\]
\[
\Rightarrow \left( {\begin{array}{*{20}c}
1 & 0 & -5 &  z^{2}+4z+5 \\
0 & 1 & \frac{25}{4}\left( z^{2}+1\right)&  -5  \\
0 & 0 & 0 &  1 \\
0 & 0 & 1 &  0\\
\end{array} } \right)
\]
\[
\Rightarrow \left( {\begin{array}{*{20}c}
1 & 0 & -5&  z^{2}+4z+5 \\
\frac{5}{4}\left( z^{2}+1\right) & 1 & 0  &  \frac{5}{4}(z^{2}+1)(z^{2}+4z+5)-5  \\
0 & 0 & 0 &  1 \\
0 & 0 & 1 &  0\\
\end{array} } \right)
\]
\[
\Rightarrow \left( {\begin{array}{*{20}c}
1 & 0 & 5 & 0 \\
\frac{5}{4}\left( z^{2}+1\right) & 1 & 0  &  \frac{5}{4}(z+1)^{4}  \\
0 & 0 & 0 &  1 \\
0 & 0 & -1 &   \frac{1}{5}\left( z^{2}+4z+5 \right)\\
\end{array} } \right).
\]
Hence,
\[
S\left( z \right)=\left( {\begin{array}{*{20}c}
1 & 0\\
\frac{5}{4}\left( z^{2}+1 \right) & 1\\
\end{array} } \right),\, D\left( z \right)=\left( {\begin{array}{*{20}c}
5 & 0\\
0 & \frac{5}{4}\left( z+1 \right)^{4}\\
\end{array} } \right),
\]
\[
T\left( z \right)=\left( {\begin{array}{*{20}c}
0 & 1\\
-1 & \frac{1}{5}\left( z^{2}+4z+5 \right)\\
\end{array} } \right),\, \tilde{D}\left( z \right)=\left( 
{\begin{array}{*{20}c}
\frac{5}{z^{2}} & 0\\
0 & \frac{5}{4z^{2}}\left( z+1 \right)^{4}\\
\end{array}}\right).
\]
According to Theorem \ref{theorem412} (i), the root functions of $Q(z)$ are among the columns of the matrix $T\left( z \right)$. We can easily verify that only 
\[
\varphi \left( z \right)=T_{2}(z)=\left( {\begin{array}{*{20}c}
1\\
\frac{1}{5}\left( z^{2}+4z+5 \right)\\
\end{array} } \right)
\]
is a root function of $Q$ at $\alpha =-1$; it satisfies
\[
Q\left( z \right)\varphi \left( z \right)=\left( {\begin{array}{*{20}c}
0\\
\frac{5}{4z^{2}}\left( z+1 \right)^{4}\\
\end{array} } \right) \, \wedge \, L\left( z \right)\varphi \left( z \right)=\left( {\begin{array}{*{20}c}
0\\
\frac{5}{4}\left( z+1 \right)^{4}\\
\end{array} } \right).
\]
When expanded in powers of $\left( z+1 \right)$ the root function $\varphi 
\left( z \right)$ is
\[
\varphi \left( z \right)=\left( {\begin{array}{*{20}c}
1\\
\raise0.7ex\hbox{$2$} \!\mathord{\left/ {\vphantom {2 
5}}\right.\kern-\nulldelimiterspace}\!\lower0.7ex\hbox{$5$}\\
\end{array} } \right)+\left( {\begin{array}{*{20}c}
0\\
\raise0.7ex\hbox{$2$} \!\mathord{\left/ {\vphantom {2 
5}}\right.\kern-\nulldelimiterspace}\!\lower0.7ex\hbox{$5$}\\
\end{array} } \right)\left( z+1 \right)+\left( {\begin{array}{*{20}c}
0\\
\raise0.7ex\hbox{$1$} \!\mathord{\left/ {\vphantom {1 
5}}\right.\kern-\nulldelimiterspace}\!\lower0.7ex\hbox{$5$}\\
\end{array} } \right)\left( z+1 \right)^{2}.
\]
This proves that the matrix polynomial $\varphi(z)$ is a root function of both $Q(z)$ and $L(z)$ at $\alpha=-1$ of order $k=4$.

Note that the expansion (\ref{eq440}) of $\varphi(z)$ contains only $l = 3$ generalized Jordan vectors $\varphi_{j}(\alpha)$ of $L(z)$ at $\alpha = -1$. This is another example showing that a root function of a polynomial or meromorphic matrix function may have an order greater than its degree.

(iii) According to Theorem \ref{theorem412} (iii), the corresponding canonical pole 
cancellation functions of $Q(z)$ at $\beta =0$ are:
\[
\psi_{1}\left( z \right)=T_{1}\left( z \right){\tilde{d}_{1}\left( z 
\right)}^{-1}=\left( {\begin{array}{*{20}c}
0\\
-\frac{z^{2}}{5}\\
\end{array} } \right),
\]
\[
\psi_{2}\left( z \right)=T_{2}\left( z \right){\tilde{d}_{2}\left( z 
\right)}^{-1}=\left( {\begin{array}{*{20}c}
1\\
\frac{1}{5}\left( z^{2}+4z+5 \right)\\
\end{array} } \right)\frac{4z^{2}}{5\left( z+1 \right)^{4}}.
\]\hfill $\square$
\section{Comparison with the logarithmic residue theorem}\label{s6}
The logarithmic residue theorem \cite[Theorem 2.1]{GLR} is a well-known 
classical result that addresses the orders and multiplicities of zeros 
and poles of operator-valued functions. It is straightforward to verify that 
this theorem applies to $Q(z)$, a meromorphic matrix function on $\bar{\mathbb{C}}$, 
provided that $\det {Q\left( z \right)}\not\equiv 0$.

If $\alpha \in \sigma_{0}\left( Q \right)$ and $s_{1}(\alpha ),\, 
\mathellipsis ,\, s_{l}\left( \alpha \right)$, are orders of the zero 
$\alpha $ corresponding to the functions $\tilde{d}_{i}(z)$ in (\ref{eq410}), respectively, then, according to  Theorem \ref{theorem412}, the total \textit{multiplicity of the zero} $\alpha$ is $N\left( \alpha \right):=s_{1}\left( \alpha \right)+\, 
\mathellipsis +\, s_{l}\left( \alpha \right)$. Similarly, if $\beta \in 
\sigma_{p}\left( Q \right)$ and $t_{1}(\beta ),\, \mathellipsis ,\, 
t_{m}\left( \beta \right)$, are orders of the pole $\beta $ 
corresponding to the functions $\tilde{d}_{i}(z)$, respectively, then the total \textit{multiplicity of the pole} $\beta $ is
$P\left( \beta \right):=t_{1}(\beta )+\, \mathellipsis +\, t_{m}(\beta )$. 
According to the logarithmic residue theorem, it holds 
\begin{equation}
\label{eq62}
\frac{1}{2\pi i}sp\oint_\Gamma {Q^{'}\left( z \right)Q^{-1}\left( z \right)} 
dz=N\left( \alpha \right)-P\left( \beta \right),
\end{equation}
where $\Gamma \subseteq \mathcal{D}\left( Q \right)$ is a simple rectifiable contour 
enclosing only $\alpha $ and $\beta $ as critical points, and ``$sp$'' denotes the trace of a matrix. 
If $\alpha \ne \beta $, and only $\alpha $ or $\beta $ lies within the 
contour $\Gamma $, then formula (\ref{eq62}) gives the total multiplicity of the 
zero $\alpha $ or of the pole $\beta $, respectively. However, this approach 
does not allow us to determine the partial zero multiplicities $s_{1}\left( 
\alpha \right),\, \mathellipsis ,s_{l}\left( \alpha \right)$, or the partial pole 
multiplicities $t_{1}(\beta ),\, \mathellipsis ,\, t_{m}\left( \beta \right)$. The 
formula (\ref{eq62}) is even less informative in the case $\alpha =\beta $. 

In the next example, we demonstrate how the method developed in this paper compares with the logarithmic residue theorem in identifying the orders and multiplicities of zeros and poles.

\begin{example}\label{example62} Find the multiplicities of the zeros and poles of the 
following function:
\[
Q\left( z \right)=\left( {\begin{array}{*{20}c}
0 & 1 & 0\\
\frac{1}{z} & 0 & 1\\
0 & 0 & z^{2}-z\\
\end{array} } \right)
\]
\end{example}
Let us first compute the integral in (\ref{eq62}) around the critical point $\alpha 
=1$, i.e., when $\Gamma $ is the circle $\left| z-1 \right|=r<1$.
\[
Q^{'}\left( z \right)=\left( {\begin{array}{*{20}c}
0 & 0 & 0\\
\frac{-1}{z^{2}} & 0 & 0\\
0 & 0 & 2z-1\\
\end{array} } \right),\, \, Q^{-1}(z)=\left( {\begin{array}{*{20}c}
0 & z & \frac{-1}{z-1}\\
1 & 0 & 0\\
0 & 0 & \frac{1}{z^{2}-z}\\
\end{array} } \right)\Rightarrow 
\]
\[
Q^{'}\left( z \right)Q^{-1}\left( z \right)=\left( {\begin{array}{*{20}c}
0 & 0 & 0\\
0 & \frac{-1}{z} & \frac{-1}{z^{2}}+\frac{1}{z-1}-\frac{1}{z}\\
0 & 0 & \frac{1}{z-1}+\frac{1}{z}\\
\end{array} } \right)
\]
\[
\Rightarrow \, \frac{1}{2\pi i}sp\oint_{\left| z-1 \right|=r} {Q^{'}\left( z 
\right)Q^{-1}\left( z \right)} dz=sp\left( {\begin{array}{*{20}c}
0 & 0 & 0\\
0 & 0 & 1\\
0 & 0 & 1\\
\end{array} } \right)=1.
\]
Since $Q\left( z \right)$ is holomorphic at $\alpha =1$, we conclude that 
the critical point $\alpha =1$ is a zero of multiplicity $s\left( 1 
\right)=1$. 

Similarly
\[
\frac{1}{2\pi i}sp\oint_{\left| z \right|=r} {Q^{'}\left( z 
\right)Q^{-1}\left( z \right)} dz=sp\left( {\begin{array}{*{20}c}
0 & 0 & 0\\
0 & -1 & -1\\
0 & 0 & 1\\
\end{array} } \right)=0.
\]
Hence, the integral provides no information about the nature 
of the point $\alpha =0$. Since $\det {Q\left( z \right)=1-z}$ 
there is no indication from the determinant either that $\alpha =0$ is either a 
zero or a pole of $Q\left( z \right)$.

Let us now apply the method presented in Section \ref{s4} to determine the zeros 
and poles of $Q\left( z \right)$, and their partial and total 
multiplicities.
\[
q\left( z \right)=z\Rightarrow L\left( z \right)=q\left( z 
\right)Q\left( z \right)=\left( {\begin{array}{*{20}c}
0 & z & 0\\
1 & 0 & z\\
0 & 0 & z^{3}-z^{2}\\
\end{array} } \right).
\]
The factorization (\ref{eq48}) is:
\[
\left( {\begin{array}{*{20}c}
0 & 1 & 0\\
1 & 0 & 0\\
0 & 0 & 1\\
\end{array} } \right)\left( {\begin{array}{*{20}c}
0 & 1 & 0\\
\frac{1}{z} & 0 & 1\\
0 & 0 & z^{2}-z\\
\end{array} } \right)\left( {\begin{array}{*{20}c}
1 & 0 & -z\\
0 & 1 & 0\\
0 & 0 & 1\\
\end{array} } \right)=\left( {\begin{array}{*{20}c}
\frac{1}{z} & 0 & 0\\
0 & 1 & 0\\
0 & 0 & z^{2}-z\\
\end{array} } \right)=\tilde{D}(z).
\]

According to Definition \ref{definition410}, $\alpha =0$ is both a pole and a zero 
of $Q$, each with multiplicity one; we have a complete description of 
the zeros and poles of $Q$, including their multiplicities.

Note that $T_{2}\left( z \right)=\left( 
{\begin{array}{*{20}c}
0\\
1\\
0\\
\end{array} } \right)$ is a root function of $L$ at $\alpha =0$, 
but it is not a root function of $Q$ at all. Such a reduction in the partial multiplicity of a zero or pole occurs only at points that are both zeros and poles.  \hfill $\square$

\section{An application of meromorphic matrix functions to differential equations}\label{s8}
In this section, we define a meromorphic matrix function $Q(z)$ associated with the nonlinear system of differential equations given by~(\ref{eq82}) below, which generalizes system~(\ref{eq22}).

We then solve system~(\ref{eq82}) using the zeros and corresponding eigenvectors of the meromorphic matrix function associated with it. All zeros~$\alpha$ of~$Q(z)$ and their corresponding eigenvectors~$\varphi(\alpha)$ can be determined using the method developed in Section~\ref{s4}, without solving the equation $\det Q(z)=0$.
\\

Let us consider the following nonlinear system of order $m \in \mathbb{N}$ consisting of 
$n\in \mathbb{N}$ differential equations with $n$ unknown functions $u_{1}\left( t 
\right),\mathellipsis ,\, u_{n}\left( t \right)$
\begin{equation}
\label{eq82}
A_{m}\left( {\begin{array}{*{20}c}
{u_{1}^{(m)}}^{-1}\\
\vdots \\
{u_{n}^{(m)}}^{-1}\\
\end{array} } \right)+A_{m-1}\left( {\begin{array}{*{20}c}
{u_{1}^{(m-1)}}^{-1}\\
\vdots \\
{u_{n}^{(m-1)}}^{-1}\\
\end{array} } \right)+\mathellipsis +A_{0}\left( {\begin{array}{*{20}c}
{u_{1}}^{-1}\\
\vdots \\
{u_{n}}^{-1}\\
\end{array} } \right)=0,
\end{equation}
where $A_{i},\, i=0,\, \mathellipsis ,\,m,\, \, A_{m}\ne 0,\, $ are 
$n\times n$ complex matrices.  Additionally, we assume $m>1$ or $n>1$ since the case $m=n=1$ reduces to a trivial case consisting of a single equation of the first order that separates variables.

We will use the following notation for vectors
\[
\varphi^{-1}:=\left( {\begin{array}{*{20}c}
\varphi_{1}^{-1}\\
\vdots \\
\varphi_{n}^{-1}\\
\end{array} } \right).\, 
\]
\begin{theorem}\label{theorem82}. Consider the system (\ref{eq82}) with $m > 1$ or $n > 1$. The function $Q$ associated with this system is defined by the expression:
\begin{equation}
\label{eq84}
Q\left( z \right)=\frac{1}{z^{m}}\left( A_{m}+A_{m-1}z+\mathellipsis +A_{0}z^{m}\right).
\end{equation}
 \begin{enumerate}[(i)]
\item The solutions of system (\ref{eq82}) are vector-valued functions of the form:
\begin{equation}
\label{eq86}
u_{j}\left( t \right)=T_{j}\left( \alpha \right)^{-1}e^{\alpha t}, \, \alpha \in \sigma_{0}\left( Q \right) , \, j \in \Omega_{0}\left( \alpha \right),
\end{equation}
where $\sigma_{0}\left( Q \right)$ and $\Omega_{0}\left( \alpha \right)$ are defined as in Section \ref{s4}, and $T_{j}(z)$ denote the columns of the matrix of elementary transformations $T(z)$ in identity (\ref{eq48}). 
\item For any linear combination 
\begin{equation}
\label{eq87}
\varphi(\alpha)= \sum\limits_{j\in \Omega_{o}(\alpha)} c_{j}T_{j}(\alpha), \, c_{j}\in \mathbb{C},
\end{equation}
we obtain a solution of system (\ref{eq82}) of the form
\begin{equation}
\label{eq88}
u\left( t \right)=\varphi\left( \alpha \right)^{-1}e^{\alpha t}.
\end{equation}
\item The number of linearly independent solutions (\ref{eq88}) of the system (\ref{eq82}) is equal to the \textbf{geometric} multiplicity of the zero $\alpha$ of the function $Q$ given by (\ref{eq84}). However, linear combinations of different solutions of the system (\ref{eq82}), including those corresponding to the same eigenvalue $\alpha \in \sigma_{0}\left( Q \right)$,
\[
\sum\limits_{j\in \Omega_{o}(\alpha)} c_{j}u_{j}(t), \, c_{j}\in \mathbb{C},
\]
are, in general, not solutions of the system (\ref{eq82}).
\end{enumerate}
\end{theorem}
\textbf{Proof}. The proof is based on the theory presented in Section \ref{s4}.
	
(i) We seek a solution of (\ref{eq82}) in the form of the vector function
\begin{equation}
\label{eq810}
u\left( t \right)=\left( {\begin{array}{*{20}c}
\varphi_{1}^{-1}\\
\vdots \\
\varphi_{n}^{-1}\\
\end{array} } \right)e^{zt}.
\end{equation}
Substituting (\ref{eq810}) into equation (\ref{eq82}), we obtain:
\[
A_{m}\left( {\begin{array}{*{20}c}
\varphi_{1}z^{-m}e^{-zt}\\
\vdots \\
\varphi_{n}z^{-m}e^{-zt}\\
\end{array} } \right)+\mathellipsis +A_{0}\left( {\begin{array}{*{20}c}
\varphi_{1}e^{-zt}\\
\vdots \\
\varphi_{n}e^{-zt}\\
\end{array} } \right)=0
\]
which simplifies to
\begin{equation}
\label{eq812}
\frac{1}{z^{m}}\left( A_{m}+A_{m-1}z+\mathellipsis +A_{0}z^{m}\right)\left( 
{\begin{array}{*{20}c}
\varphi_{1}\\
\vdots \\
\varphi_{n}\\
\end{array} } \right)=0.
\end{equation}
This equation leads to definition (\ref{eq84}) of the function $Q(z)$. 

According to Theorem \ref{theorem412} (i), for $i \in \Omega_{0}(\alpha) $, the vector $\varphi(\alpha) = T_{i}(\alpha)$ satisfies the algebraic system (\ref{eq812}) at $z=\alpha$, since it is an eigenvector of the matrix function $Q(z)$ defined by (\ref{eq84}). This is equivalent to stating that the vector-valued functions given in (\ref{eq86}) satisfy system (\ref{eq82}), thereby completing the proof of part (i).    

(ii) The vector $\varphi(\alpha)$, defined by (\ref{eq87}), is a linear combination of eigenvectors $T_{j}(\alpha), \, j\in \Omega_{0}(\alpha)$, and is therefore itself an such eigenvector. Thus, the vector $\varphi(\alpha)$ satisfies equation (\ref{eq812}) for $z=\alpha$. Consequently, the function $u\left( t \right)=\varphi\left( \alpha \right)^{-1}e^{\alpha t}$ is a solution of system (\ref{eq82}).

(iii) By definition, the number of linearly independent eigenvectors in (\ref{eq87}) equals the geometric multiplicity of the zero $\alpha$ of $Q$, which is also the number of independent solutions of (\ref{eq88}).

Here, \textit{independent solutions} means that the vectors $\varphi_{i}(\alpha)$ form linearly independent combinations of the type (\ref{eq87}). In fact, a linear combination of the solutions (\ref{eq86}) of system (\ref{eq82}) is not necessarily another solution of (\ref{eq82}) when $m > 1$ or $n > 1$, even if the solutions correspond to the same eigenvalue $\alpha$. This occurs because
\[
\alpha \neq \beta \Rightarrow T_{i}\left( \alpha \right)^{-1}e^{-\alpha t}+T_{i}\left( \beta \right)^{-1}e^{-\beta t}\ne \left( 
T_{i}\left( \alpha \right)e^{\alpha t}+T_{i}\left( \beta \right)e^{\beta t} \right)^{-1}.
\]
and
\[
i \neq j \Rightarrow T_{i}\left( \alpha \right)^{-1}+T_{j}\left( \alpha \right)^{-1}\ne \left( 
T_{i}\left( \alpha \right)+T_{j}\left( \alpha \right) \right)^{-1}.
\]
\hfill $\square$

\begin{example}\label{example84} Solve the nonlinear system 
\[
\frac{1}{u_{1}^{'}}+\frac{1}{u_{2}^{''}}=0,
\]
\begin{equation}
\label{eq814}
\frac{1}{u_{1}}+\frac{4}{u_{1}^{'}}+\frac{4}{u_{1}^{''}}+\frac{1}{u_{2}}+\frac{4}{u_{2}^{'}}+\frac{4}{u_{2}^{''}}=0.
\end{equation}
\end{example}
We will follow the theory and notations developed above. We seek the 
solution in the form (\ref{eq810}), which leads us to the corresponding matrix 
function (\ref{eq84}): 
\[
Q\left( z \right)=\left( {\begin{array}{*{20}c}
\frac{1}{z} & \frac{1}{z^{2}}\\
1+\frac{4}{z}+\frac{4}{z^{2}} & 1+\frac{4}{z}+\frac{4}{z^{2}}\\
\end{array} } \right)\Rightarrow 
\]
\[
L\left( z \right)=\left( {\begin{array}{*{20}c}
z & 1\\
z^{2}+4z+4 & z^{2}+4z+4\\
\end{array} } \right)\Rightarrow 
\]
\[
\Rightarrow \left( {\begin{array}{*{20}c}
I_{11} & L(z)\\
0 & I_{22}\\
\end{array} } \right)=\left( {\begin{array}{*{20}c}
1 & 0 & z & 1\\
0 & 1 & z^{2}+4z+4 & z^{2}+4z+4\\
0 & 0 & 1 & 0\\
0 & 0 & 0 & 1\\
\end{array} } \right)
\]
\[
\Rightarrow \left( {\begin{array}{*{20}c}
1 & 0 & 0 & 1\\
0 & 1 & ( 1-z)(z^{2}+4z+4) & z^{2}+4z+4\\
0 & 0 & 1 & 0\\
0 & 0 & -z & 1\\
\end{array} } \right)
\]
\[
\Rightarrow \left( {\begin{array}{*{20}c}
1 & 0 & 0 & 1\\
-z^{2}-4z-4 & 1 & ( 1-z)(z^{2}+4z+4) & 0\\
0 & 0 & 1 & 0\\
0 & 0 & -z & 1\\
\end{array} } \right)
\]
\[
\Rightarrow \left( {\begin{array}{*{20}c}
1 & 0 & 1 & 0 \\
-(z+2)^{2} & 1 & 0 & ( 1-z)(z+2)^{2}\\
0 & 0 & 0 & 1\\
0 & 0 & 1 & -z\\
\end{array} } \right)
\]
\[
\Rightarrow S\left( z \right)=\left( {\begin{array}{*{20}c}
1 & 0\\
-(z+2)^{2} & 1\\
\end{array} } \right),\, D\left( z \right)=\left( {\begin{array}{*{20}c}
1 & 0\\
0 & \left( 1-z \right)\left( z+2 \right)^{2}\\
\end{array} } \right),
\]
\[
T\left( z \right)=\left( {\begin{array}{*{20}c}
0 & 1\\
1 & -z\\
\end{array} } \right),\, \tilde{D}\left( z \right)=\left( 
{\begin{array}{*{20}c}
\frac{1}{z^{2}} & 0\\
0 & \frac{\left( 1-z \right)\left( z+2 \right)^{2}}{z^{2}}\\
\end{array} } \right).
\]
Obviously, $T_{1}=\left( {\begin{array}{*{20}c}
0\\
1\\
\end{array} } \right)$ is not a root function of $Q(z)$, while $T_{2}\left( 
z \right)=\left( {\begin{array}{*{20}c}
1\\
-z\\
\end{array} } \right)$ is a root function of $Q(z)$ at both zeros $\alpha 
=1$ and $\alpha =-2$. Thus, we have two eigenvectors and two corresponding 
solutions $u^{1}\left( z \right)$ and $u^{2}\left( z \right)$ of the system 
(\ref{eq814}), i.e.
\[
T_{2}\left( 1 \right)=\left( {\begin{array}{*{20}c}
1\\
-1\\
\end{array} } \right)\Rightarrow u^{1}\left( z \right)=\left( 
{\begin{array}{*{20}c}
1\\
-1\\
\end{array} } \right)e^{t},
\]
\[
T_{2}\left( -2 \right)=\left( {\begin{array}{*{20}c}
1\\
2\\
\end{array} } \right)\Rightarrow u^{2}\left( z \right)=\left( 
{\begin{array}{*{20}c}
1\\
\raise0.7ex\hbox{$1$} \!\mathord{\left/ {\vphantom {1 
2}}\right.\kern-\nulldelimiterspace}\!\lower0.7ex\hbox{$2$}\\
\end{array} } \right)e^{-2t}.
\]\hfill $\square$
\section{A new operator realization of a function $Q\in \mathcal{N}_{\kappa}(\mathcal{H})$}\label{s10}
\textbf{5.1.} Let $\kappa \in \mathbb{N}\cup \left\{ 0 \right\}$ and let $\left( \mathcal{K},\, \left[ .,. 
\right] \right)$ be a Krein space. Recall, it is a complex vector space equipped 
with a Hermitian sesquilinear form $\left[ \mathrm{.,.} \right]$ such that 
the following orthogonal decomposition exists:
\[
\mathcal{K}=\mathcal{K}_{+}+\mathcal{K}_{-},
\]
where $\left( \mathcal{K}_{+},[.,.] \right)$ and $\left( 
\mathcal{K}_{-},-[.,.] \right)$ are Hilbert spaces which are 
mutually orthogonal with respect to the form $[.,.]$. If the scalar 
product $\left[ .,. \right]$ has $\kappa (<\infty )$ negative squares, then 
it is called a Pontryagin space of (negative) index $\kappa $. A bounded 
operator $A$ is called self-adjoint in the Pontryagin space $\left( 
\mathcal{K},\left[ .,. \right] \right)$ if it satisfies 
\[
\left[ Ax,y \right]=\left[ x,Ay \right],\, \forall x,\, y\in \mathcal{K}
\]
Other definitions and properties related to Krein and Pontryagin spaces can be found e.g. in \cite{Bog,IKL}.

The first statement of \cite[Lemma 3]{B2} we use in this section. For convenience if the reader we repat it here as:
\begin{lemma}\label{lemma101} A function $Q\in \mathcal{N}_{\kappa }(\mathcal{H})$ is holomorphic at $\infty$ with $Q \left( \infty \right):=0 $ if and only if $Q(z)$ has minimal representation 
\[
Q\left( z \right)=\Gamma^{+}\left( A-z \right)^{-1}\Gamma, \thinspace z\in \mathcal{D}(Q), 
\]
with a bounded self-adjoint operator $A$ in a minimal Pontryagin space $\mathcal{K}$, a bounded operator $\Gamma:\mathcal{H}\to \mathcal{K}$ and its adjoint operator $\Gamma^{+}:\, \mathcal{K}\to \mathcal{H}$.
\end{lemma} 
\begin{lemma}\label{lemma102}. A function $Q\in \mathcal{N}_{\kappa }(\mathcal{H})$ is holomorphic at 
$\infty $ if and only if $Q(z)$ admits the representation 
\begin{equation}
\label{eq101}
Q\left( z \right)=S+\Gamma^{+}\left( A-z \right)^{-1}\Gamma ,\, S^{\ast }=S\in \mathcal{L}\left( \mathcal{H} \right),
\end{equation} 
with a self-adjoint bounded operator $A$ in the minimal Pontryagin space $\mathcal{K}$, 
that is,
\begin{equation}
\label{eq102}
\mathcal{K}=c.l.s.\left\{ \left( A-z \right)^{-1}\Gamma h:z\in \mathcal{D}(Q),h\in \mathcal{H}\right\},
\end{equation}
where $\Gamma :\, \mathcal{H}\to \mathcal{K}$ is a bounded operator and $\Gamma^{+}:\, \mathcal{K}\to \mathcal{H}$ is its adjoint operator. 
\end{lemma}
\textbf{Proof.} If $Q\in \mathcal{N}_{\kappa }(\mathcal{H})$ is holomorphic at $\infty$, then the limit
\[
Q(\infty):=S=\lim\limits_{z\to \infty}{Q\left( z \right)},
\]
exists. This means that the function 
$
Q_{1}\left( z \right):=Q\left( z \right)-S
$
is holomorphic at $\infty$ with $Q_{1}\left( \infty \right):=0$. According to Lemma \ref{lemma101}, there exists a self-adjoint bounded operator $A$ in a minimal Pontryagin space $\mathcal{K}$ 
with negative index $\kappa $, such that the function $Q_{1}\left( z \right)$ admits the minimal representation 
\[
Q_{1}\left( z \right):=Q\left( z \right)-S=\Gamma^{+}\left( A-z \right)^{-1}\Gamma.
\]
This proves (\ref{eq101}). The relation (\ref{eq102}) follows from Lemma \ref{lemma101}. 

The converse statement also follows easily form Lemma \ref{lemma101}. \hfill $\square$ 

If a function $Q\in \mathcal{N}_{\kappa }(\mathcal{H})$ has a pole at $\infty $, then the Krein-Langer representation is represented by  a linear relation $A$ rather than by a bounded operator. We will now prove that $Q(z)$ can be represented in a slightly modified form using a bounded operator $\tilde{A}$ rather than the linear relation $A$. 

\begin{theorem}\label{theorem104}.
\begin{enumerate}[(i)]
\item If the function $Q\in \mathcal{N}_{\kappa }(\mathcal{H})$ has a pole at $\infty $ of order $m\in \mathbb{N}$, then it admits the factorization 
\begin{equation}
\label{eq108}
Q\left( z \right)=\left( z-\beta \right)^{m}\tilde{Q}\left( z \right),
\end{equation}
where $\beta \in \mathbb{R}$ is a regular point of $Q$, the function $\tilde{Q}\in \mathcal{N}_{\tilde{\kappa}}(\mathcal{H})$ has a pole at $\beta$ of order $m$ and admits the representation (\ref{eq101}), that is
\begin{equation}
\label{eq109}
\tilde{Q}\left( z \right)=S+\tilde{\Gamma}^{+}\left( A-z \right)^{-1}\tilde{\Gamma} ,\, S^{\ast }=S\in \mathcal{L}\left( \mathcal{H} \right),
\end{equation}
with a self-adjoint bounded operator $\tilde{A}$ in the minimal Pontryagin space 
\[
\tilde{\mathcal{K}}=c.l.s.\left\{ \left( \tilde{A}-z \right)^{-1}\tilde{\Gamma }h:z\in \mathcal{D}(\tilde{Q}),h\in (\mathcal{H}) \right\}
\]
and linear operators $\tilde{\Gamma }:\mathcal{H}\to \tilde{\mathcal{K}}$ and $\tilde{\Gamma}^{+}:\tilde{\mathcal{K}}\to \mathcal{H}$. 

Conversely, if the function $\tilde{Q} \in \mathcal{N}_{\tilde{\kappa}}(\mathcal{H})$ with a pole of order $m$ at $\beta \in \mathbb{R}$ admits the representation (\ref{eq109}), then the function $Q(z)$ given by (\ref{eq108}) has a pole at $\infty$ of order $m$ and satisfies $Q\in \mathcal{N}_{\kappa }(\mathcal{H})$, for some $\kappa \in \mathbb{N}\cup \left\lbrace 0\right\rbrace $.
\item If $A$ is the linear relation representing the function $Q$ in the Krein-Langer-type representation, and if $\kappa_{\infty }$ denotes the negative index of $A(0)$, while $\kappa_{\beta }$ denotes the negative index of the root space of $\tilde{A}$ at $\beta $, then
\begin{equation}
\label{eq1010}
\tilde{\kappa}=\kappa +\kappa_{\beta }-\kappa_{\infty }.
\end{equation}
\item If $d_{\infty}:= \dim \tilde{A}(0)$ and $d_{\beta}$ denotes the dimension of the root space of $\tilde{A}$ corresponding to the eigenvalue $\beta$, then in the general case 
\[
\kappa_{\infty }\neq \kappa_{\beta } \, \wedge \, d_{\infty} \neq d_{\beta}.
\]
\end{enumerate}
\end{theorem}
\textbf{Proof.} 

(i) Assume that $\infty $ is a pole of $Q\in \mathcal{N}_{\kappa }(\mathcal{H})$ of order $m\in \mathbb{N}$. Then for a regular point $\beta \in \mathbb{R}$ of $Q$ we have 
\[
\lim\limits_{z\to \infty}\frac{Q(z)}{\left( z-\beta \right)^{m}}=:S\neq 0.
\]

By definition, $\beta $ being a regular point of $Q$ means that $Q$ is holomorphic at $\beta$ and $Q(\beta)$ is boundedly invertible. Since we have chosen $\beta \in \mathbb{R}$, the function 
\[
\tilde{Q}\left( z \right)=\frac{Q\left( z \right)}{\left( z-\beta \right)^{m}}
\]
belongs to the class $\mathcal{N}_{\tilde{\kappa}}(\mathcal{H})$ for some $\tilde{\kappa}\in \mathbb{N}\cup \left\{ 0 \right\}$. Since $Q^{-1}(\beta)\ne 0$ is a bounded operator, $\beta$ is a pole of $\tilde{Q}$ of exact order $m$. Clearly, $\tilde{Q}$ satisfies all assumptions of Lemma \ref{lemma102}. Therefore, we have
\[
\tilde{Q}\left( z \right)=\frac{Q\left( z \right)}{\left( z-\beta \right)^{m}}=S+\tilde{\Gamma }^{+}\left( \tilde{A}-z \right)^{-1}\tilde{\Gamma },\, \, \, S^{\ast }=S\in L\left( \mathcal{H}\right),
\]
\[
\tilde{\mathcal{K}}=c.l.s.\left\{ \left( \tilde{A}-z \right)^{-1}\tilde{\Gamma }h:z\in 
\mathcal{D}(\tilde{Q}),h\in  \mathcal{H} \right\}.
\]
This completes the proof of part (i) in one direction. 

Conversely, since $\tilde{Q}\left( z \right)$ has representation (\ref{eq109}), according to Lemma \ref{lemma102}, the function $\tilde{Q}\left( z \right)$ is holomorphis at $\infty$ and $Q(\infty)=S$. Now the converse statement follows form (\ref{eq108}).

(ii) Let $\mathcal{K}$ now denote the Pontryagin space associated with the Krein–Langer representation of $Q$; see \cite[Proposition 1]{Lu3} or \cite[Theorem 1.1]{B3}. To prove identity (\ref{eq1010}), we recall the decomposition of the Pontryagin space $\mathcal{K}$, given in \cite[Proposition 5.1]{B3}
\begin{equation}
\label{eq1012}
\mathcal{K}=\mathcal{K}_{0}\left[ + \right]\mathcal{K}_{1}\left[ + \right]\mathellipsis \left[ + \right]\mathcal{K}_{r}\left[ + \right]\mathcal{K}_{r+1},
\end{equation}
\textit{where $r\in \mathbb{N}$ is the number of independent non-degenerate Jordan chains of the representing relation $A$ of $Q$ at a generalized pole $\alpha ;$ the sub-spaces $\mathcal{K}_{i},\, i=0,\, 1,\, \mathellipsis ,\, r,\, r+1$, are A-invariant Pontryagin subspaces with respective indices $\kappa_{i}.$ Moreover, $\kappa_{0}=0$ and $\kappa =\sum\limits_{i=1}^{r+1} \kappa_{i} $. For each $i=1,\, 2,\, \mathellipsis ,\, r$, the subspace $\mathcal{K}_{i}$ is a linear span of the corresponding maximal non-degenerate Jordan chain ${ x}_{0}^{i},\, \mathellipsis ,\, x_{l_{i}-1}^{i}$. All positive eigenvectors form the subspace $\mathcal{K}_{0}$. All degenerate chains of $A$ at $\alpha $ are contained in $\mathcal{K}_{r+1}$.} 

Obviously, the eigenvectors and Jordan chains corresponding to all other eigenvalues of the representing relation $A$ lie in the Pontryagin space $\mathcal{K}_{r+1}$. Therefore, we can apply the decomposition (\ref{eq1012}) on $\mathcal{K}_{r+1}$ using the next eigenvalue $\lambda \neq \alpha$. Because $A$ has only finitely many eigenvalues, iterating this procedure over all of them ultimately yields the decomposition
\begin{equation}
\label{eq1014}
\mathcal{K}=\mathcal{K}_{0}\left[ + \right]\mathcal{K}_{1}\left[ + \right]\mathellipsis \left[ + \right]\mathcal{K}_{s}\left[ + \right]\mathcal{K}_{s+1},
\end{equation}
where $s\in \mathbb{N}$ is the number of \textbf{all} linearly independent non-positive eigenvectors of $A$ with non-degenerate Jordan chains, and $\mathcal{K}_{s+1}$ contains all degenerate chains of $A$. Hence, $\kappa =\sum\limits_{i=1}^{s+1} \kappa_{i} $.

Since $\infty$ is by assumption a pole of $Q$, the subspace $A(0)$ is non-degenerate. Thus, we may take $\mathcal{K}_{s}=A(0)$, and consequently, $\kappa_{s }:=\kappa_{\infty }$. 

From the representation (\ref{eq108}), we observe that the critical and regular points of $Q$ and $\tilde{Q}$ coincide, except for $\beta$ and $\infty$. This implies that the relation $A$ and operator $\tilde{A}$ share the same eigenvalues and Jordan chains, with the same structural properties, except at $\infty$ and $\beta$. Therefore, if the subspaces $\tilde{\mathcal{K}}{j}$, $j = 1, \ldots, s+1$, play the same role for $\tilde{A}$ as the subspaces $\mathcal{K}{j}$, $j = 1, \ldots, s+1$, do for $A$, we can identify the subspaces as follows:
\[
\tilde{\mathcal{K}}_{0}=\mathcal{K}_{0},\, \mathellipsis, \, \tilde{\mathcal{K}}_{s-1}=\mathcal{K}_{s-1},\, \tilde{\mathcal{K}}_{s+1}=\mathcal{K}_{s+1},
\]
Since $\beta$ is a pole of $\tilde{Q}$, the subspace $\tilde{\mathcal{K}}_{\beta}$ is non-degenerate. Thus, we may set $\tilde{\mathcal{K}}_{s}=\tilde{\mathcal{K}}_{\beta}$. Then, the decomposition (\ref{eq1012}) of $\tilde{K}$ becomes:
\[
\tilde{\mathcal{K}}=\mathcal{K}_{0}\left[ + \right]\mathellipsis \left[ + \right]\mathcal{K}_{s-1}\left[ + \right]\tilde{\mathcal{K}}_{\beta}\left[ + \right]\mathcal{K}_{s+1}.
\]
Then, using this decomposition and the decomposition (\ref{eq1014}), we obtain:
\[
\left( \tilde{\kappa}=\sum\limits_{i=1}^{s-1} \kappa_{i} +\kappa_{\beta }+\kappa_{s+1}\, 
\, \wedge \,\, \kappa =\sum\limits_{i=1}^{s-1} \kappa_{i} +\kappa_{\infty}+\kappa_{s+1} 
\right)\Rightarrow \, \tilde{\kappa}-\kappa_{\beta }=\kappa -\kappa_{\infty }.
\]
This proves (\ref{eq1010}). 

(iii) Both claims will be demonstrated in Example~\ref{example105}. \hfill $\square$ 
\\

Since the order of the pole $\beta$ of $\tilde{Q}$ is the same as the order of the pole at $\infty$ of $Q$, the question arises whether the dimensions or negative indices of the corresponding root subspaces also coincide. The following example shows that neither of these properties holds in the general case. Consequently, the state spaces $\mathcal{K}$ and $\mathcal{\tilde{K}}$ are structurally different and have different negative indeces $\kappa$ and $\tilde{\kappa}$.

\begin{example}\label{example105} Using the matrix polynomial $L\left( z \right)=\left( {\begin{array}{*{20}c}
0 & z & 0\\
z & z^{2} & 0 \\
0 & 0 & z \\
\end{array} } \right) $ and $\beta=1$ show that, in the general case, $\kappa_{\beta} \neq \kappa_{\infty}$ and that the root sub-space of $\tilde{A}$ at $\beta$ and  $A(0)$ have different dimensions.
\end{example}
Here, the polynomial $L$ plays the role of $Q$ in Theorem \ref{theorem104}. It has a pole at $\infty$ of order $m=2$. By definition, investigating the pole at $\infty$ of $L(z)$ is equivalent to investigating the pole at $0$ of the function 
\[
L_{0}(\zeta):=-L\left( \frac{1}{\zeta}\right) =\left( {\begin{array}{*{20}c}
0 & -\frac{1}{\zeta} & 0\\
 -\frac{1}{\zeta} & -\frac{1}{\zeta^{2}}  & 0 \\
0 & 0 & -\frac{1}{\zeta} \\
\end{array} } \right).
\]
It is easy to see that the Krein–Langer representation of $L_{0}(\zeta)$ has the form 
\[
L_{0}(\zeta)=\Gamma^{+}\left( A-\zeta\right)^{-1}\Gamma, 
\]
where
\[
A=\left( {\begin{array}{*{20}c}
0 & 1 & 0\\
0 & 0 & 0 \\
0 & 0 & 0 \\
\end{array} } \right) ,\, \Gamma=I_{3}=\left( {\begin{array}{*{20}c}
1 & 0 & 0\\
0 & 1 & 0 \\
0 & 0 & 1 \\
\end{array} } \right), \, \Gamma^{+}=J=\left( {\begin{array}{*{20}c}
0 & 1 & 0\\
1 & 0 & 0 \\
0 & 0 & 1 \\
\end{array} } \right).
\]
Thus, the representing operator of $L_{0}(\zeta)$ has two Jordan chains at $0$, of maximal lengths $k_{1}=2$ and $k_{1}=2$. Hence, $d_{\infty} = 3$ and $\kappa_{\infty} = 1$. 

For $\beta = 1$, (\ref{eq108}) gives
\[
\tilde{L}\left( z \right):=\frac{L(z)}{ (z-1)^{2}}=\left( {\begin{array}{*{20}c}
0 & \frac{z}{ (z-1)^{2}} & 0\\
\frac{z}{ (z-1)^{2}} & \frac{z^{2}}{ (z-1)^{2}} & 0 \\
0 & 0 & \frac{z}{ (z-1)^{2}} \\
\end{array} } \right).
\]
In the sequel, $\tilde{L}(z)$ takes the role of $Q$ in Theorem \ref{theorem48}. Applying the diagonalization process of $\tilde{L}(z)$ described in Theorem \ref{theorem48}, we obtain the expression (\ref{eq48})
\[
\tilde{D}\left( z \right)=S(z)\tilde{L}\left( z \right)T(z)=\left( {\begin{array}{*{20}c}
1 & 0 & 0\\
0 & 1 & 0 \\
0 & 0 & 1 \\
\end{array} } \right)\left( {\begin{array}{*{20}c}
0 & \frac{z}{ (z-1)^{2}} & 0\\
\frac{z}{ (z-1)^{2}} & \frac{z^{2}}{ (z-1)^{2}} & 0 \\
0 & 0 & \frac{z}{ (z-1)^{2}} \\
\end{array} } \right) \left( {\begin{array}{*{20}c}
-z & 1 & 0\\
1 & 0 & 0 \\
0 & 0 & 1 \\
\end{array} } \right).
\]
This means,
\[
\tilde{D}\left( z \right)=\left( {\begin{array}{*{20}c}
\frac{z}{ (z-1)^{2}} & 0 & 0\\
0 & \frac{z}{ (z-1)^{2}} & 0 \\
0 & 0 & \frac{z}{ (z-1)^{2}} \\
\end{array} } \right); \, S(z)=I_{3}; \, T(z)=\left( {\begin{array}{*{20}c}
-z & 1 & 0\\
1 & 0 & 0 \\
0 & 0 & 1 \\
\end{array} } \right);
\]
\[
\tilde{D}^{-1}\left( z \right)=\left( {\begin{array}{*{20}c}
\frac{(z-1)^{2}}{ z} & 0 & 0\\
0 & \frac{(z-1)^{2}}{ z} & 0 \\
0 & 0 & \frac{(z-1)^{2}}{ z} \\
\end{array} } \right).
\]
According to Theorem \ref{theorem412} (iii) the pole cancelation functions of $\tilde{L}\left( z \right)$ at $\beta=1$ are given by the columns of the following matrix 
\[
\Psi \left( z \right):=\left[ \psi_{1}(z)\,\,\psi_{2}(z)\,\,\psi_{3}(z)\right]=  T(z)\tilde{D}\left( z \right)^{-1}=\left( {\begin{array}{*{20}c}
-(z-1)^{2} & \frac{(z-1)^{2}}{ z} & 0\\
\frac{(z-1)^{2}}{ z} & 0 & 0 \\
0 & 0 & \frac{(z-1)^{2}}{ z} \\
\end{array} } \right).
\]
Hence, $\tilde{L}(z)$ has three canonical pole cancellation functions of second order at $\beta =1$. Consequently, the representing operator of $\tilde{L}(z)$ has three canonical Jordan chains, implying $d_{\beta} = 6$ and $\kappa_{\beta} = 3 $. 

This proves Theorem \ref{theorem104} (iii). \hfill $\square$ 
\\
\\
\textbf{5.2.} As previously mentioned, obtaining a Krein–Langer-type representation for a function $Q$ with a pole at $\infty$, where $Q$ is represented by a linear relation, is generally not feasible. Although the representation (\ref{eq108}) is not a Krein–Langer representation, it allows us to work with a bounded operator $\tilde{A}$ in place of the (unbounded) linear relation $A$. In the following example, we illustrate how to apply the root and pole cancellation functions to construct the operator $\tilde{A}$ and the representation (\ref{eq108}) of a given rational function $Q \in \mathcal{N}_{\kappa}^{n\times n}$ with a pole at $\infty$.

\begin{example}\label{example106} For the function
\[
Q\left( z \right)=\left( {\begin{array}{*{20}c}
\frac{1}{z} & z-1\\
z-1 & 0\\
\end{array} } \right).
\]
find the operator representation of the form (\ref{eq108}). 
\end{example}
Obviously, $z=\infty $ is a pole of $Q\left( z \right)$ of order $m=1$. We can choose, for example, $\beta =-1$. Then we define
\[
\tilde{Q}\left( z \right):=\frac{Q\left( z \right)}{z+1}.
\]
We have
\[
S:=\lim\limits_{z\to \infty}{\tilde{Q}\left( z \right)=}\left( {\begin{array}{*{20}c}
\frac{1}{z\left( z+1 \right)} & \frac{z-1}{ z+1}\\
\frac{z-1}{ z+1 } & 0\\
\end{array} } \right)=\left( {\begin{array}{*{20}c}
0 & 1\\
1 & 0\\
\end{array} } \right)\,
\wedge \,
L\left( z \right)=\left( {\begin{array}{*{20}c}
1 & z^{2}-z\\
z^{2}-z & 0\\
\end{array} } \right).
\]
In this example, the identity (\ref{eq48}) takes the form:
\[
S\left( z \right)\tilde{Q}\left( z \right)T\left( z \right)=\left( {\begin{array}{*{20}c}
1 & 0\\
-\left( z^{2}-z \right) & 1\\
\end{array} } \right)\left( {\begin{array}{*{20}c}
\frac{1}{z\left( z+1 \right)} & \frac{z-1}{ z+1 }\\
\frac{z-1}{ z+1 } & 0\\
\end{array} } \right)\left( {\begin{array}{*{20}c}
1 & -\left( z^{2}-z \right)\\
0 & 1\\
\end{array} } \right)=
\]
\[
=\left( {\begin{array}{*{20}c}
\frac{1}{z\left( z+1 \right)} & 0\\
0 & -\frac{z\left( z-1 \right)^{2}}{z+1}\\
\end{array} } \right)=:\tilde{D}\left( z \right).
\]
Hence,
\[
S^{-1}\left( z \right)=\left( {\begin{array}{*{20}c}
1 & 0\\
z\left( z-1 \right) & 1\\
\end{array} } \right),\, \, T\left( z \right)=\left( {\begin{array}{*{20}c}
1 & -z\left( z-1 \right)\\
0 & 1\\
\end{array} } \right).
\]
Then, 
\[
T_{2}\left( z \right)=\left( {\begin{array}{*{20}c}
-z\left( z-1 \right)\\
1\\
\end{array} } \right)
\]
is a root function of $\tilde{Q}(z)$ of second order at $\alpha=1$, and of first order at $\alpha=0$. 

According to Theorem \ref{theorem412} (iii), we can find all pole cancellation functions of $\tilde{Q}$ as the columns of the following matrix:
\[
T\left( z \right)\tilde{D}^{-1}\left( z \right)=\left( 
{\begin{array}{*{20}c}
1 & -\left( z^{2}-z \right)\\
0 & 1\\
\end{array} } \right)\left( {\begin{array}{*{20}c}
z\left( z+1 \right) & 0\\
0 & -\frac{z+1}{z\left( z-1 \right)^{2}}\\
\end{array} } \right).
\]
Hence, we obtain 
\[
\left[ \psi_{1}\left( z \right)\,\,\psi_{2}\left( z \right)\, 
\right]:=\left( {\begin{array}{*{20}c}
z\left( z+1 \right) & \frac{z+1}{z-1}\\
0 & -\frac{z+1}{z\left( z-1 \right)^{2}}\\
\end{array} } \right).
\]
The function $\psi_{1}\left( z \right)$ is a pole cancellation function of $Q$ for the poles $\beta =0$ and $\beta =-1$, each of the first order. The function $\psi_{2}\left( z \right)$ is a pole cancellation function for the pole $\beta =-1$, of the first order.

According to Theorem \ref{theorem412} (iv), the functions
\[
\hat{\varphi }_{1}\left( z \right)=\hat{S}_{1}\left( z \right)=\left( 
{\begin{array}{*{20}c}
1\\
z\left( z-1 \right)\\
\end{array} } \right)\, \wedge \, \hat{\varphi }_{2}\left( z 
\right)=\hat{S}_{2}\left( z \right)=\left( {\begin{array}{*{20}c}
0\\
1\\
\end{array} } \right)\, .
\]
are the corresponding pole functions of $\tilde{Q}$. 

Using the notation of \cite[Corollary 3.10]{B5}, we have $\hat{Q}=-Q^{-1}$. Therefore, 
we apply that corollary with 
\[
\Psi_{1}\left( z \right)=\left( {\begin{array}{*{20}c}
-z\left( z+1 \right)\\
0\\
\end{array} } \right),\, \Psi_{2}\left( z \right)=\left( 
{\begin{array}{*{20}c}
-\frac{z+1}{z-1}\\
\frac{z+1}{z\left( z-1 \right)^{2}}\\
\end{array} } \right).
\]
We have:
\[
\lim\limits_{z\to \beta}\frac{\left( \Psi_{1}\left( z \right),\hat{\varphi }_{1}\left( z 
\right) \right)}{\left( z-\beta \right)^{t_{\beta }}}=\lim\limits_{z\to 0}\frac{\left( \left( 
{\begin{array}{*{20}c}
-z\left( z+1 \right)\\
0\\
\end{array} } \right),\left( {\begin{array}{*{20}c}
1\\
z\left( z-1 \right)\\
\end{array} } \right) \right)}{z}=-1<0.
\]
\[
\lim\limits_{z\to \beta}\frac{\left( \Psi_{1}\left( z \right),\hat{\varphi }_{1}\left( z \right) \right)}{\left( z-\beta \right)^{t_{\beta }}}=\lim\limits_{z\to -1}\frac{\left( \left( 
{\begin{array}{*{20}c}
-z\left( z+1 \right)\\
0\\
\end{array} } \right),\left( {\begin{array}{*{20}c}
1\\
z\left( z-1 \right)\\
\end{array} } \right) \right)}{z+1}=1>0.
\]
\[
\lim\limits_{z\to \beta}\frac{\left( \Psi_{2}\left( z \right),\hat{\varphi }_{2}\left( z \right) \right)}{\left( z-\beta \right)^{t_{\beta}}}=\lim\limits_{z\to -1}\frac{\left( \left( 
{\begin{array}{*{20}c}
-\frac{z+1}{z-1}\\
\frac{z+1}{z\left( z-1 \right)^{2}}\\
\end{array} } \right),\left( {\begin{array}{*{20}c}
0\\
1\\
\end{array} } \right) \right)}{z+1}=-\frac{1}{4}<0.
\]
According to \cite[Theorem 3.8]{B5}, the representing operator $\tilde{A}$ has three eigenvectors with no Jordan chains of higher order. Consequently, the Pontryagin space $\left( \tilde{\mathcal{K}},\left[ .,. \right] \right)$ is three-dimensional. 

The signature of the above limits determine the fundamental symmetry $\tilde{J}$ according to \cite[Lemma 3.7]{B5}. This yields the following form for the operator $\tilde{A}$ and the fundamental symmetry $\tilde{J}$:
\[
\tilde{A}=\left( {\begin{array}{*{20}c}
0 & 0 & 0\\
0 & -1 & 0\\
0 & 0 & -1\\
\end{array} } \right)\, \wedge \, \, \tilde{J}=\left( {\begin{array}{*{20}c}
-1 & 0 & 0\\
0 & 1 & 0\\
0 & 0 & -1\\
\end{array} } \right)
\]
\[
\Rightarrow \, \left( \tilde{A}-z \right)^{-1}=\left( {\begin{array}{*{20}c}
-\frac{1}{z} & 0 & 0\\
0 & \frac{-1}{z+1} & 0\\
0 & 0 & \frac{-1}{z+1}\\
\end{array} } \right).
\]
We will now solve the following equation for $\tilde{\Gamma} $:
\begin{equation}
\label{eq104}
\tilde{Q}\left( z \right)-S=\tilde{\Gamma}^{+}\left( \tilde{A}-z \right)^{-1}\tilde{\Gamma}. 
\end{equation}
We have
\[
\tilde{Q}\left( z \right)-S=\left( {\begin{array}{*{20}c}
\frac{1}{z\left( z+1 \right)} & \frac{z-1}{\left( z+1 \right)}\\
\frac{z-1}{\left( z+1 \right)} & 0\\
\end{array} } \right)-\left( {\begin{array}{*{20}c}
0 & 1\\
1 & 0\\
\end{array} } \right)=\left( {\begin{array}{*{20}c}
\frac{1}{z\left( z+1 \right)} & \frac{-2}{z+1}\\
\frac{-2}{z+1} & 0\\
\end{array} } \right).
\]
Since the Pontryagin space is 3-dimensional and $Q\left( z \right):\mathbb{C}^{2}\to 
\mathbb{C}^{2}$, we seek $\tilde{\Gamma} $ as a $3\times 2$ matrix:
\[
\tilde{\Gamma} =\left( {\begin{array}{*{20}c}
\gamma_{11} & \gamma_{12}\\
\gamma_{21} & \gamma_{22}\\
\gamma_{31} & \gamma_{32}\\
\end{array} } \right)\Rightarrow \tilde{\Gamma}^{+}=\tilde{\Gamma}^{\ast }\tilde{J}=\left( 
{\begin{array}{*{20}c}
-\bar{\gamma_{11}} & \bar{\gamma_{21}} & -\bar{\gamma_{31}}\\
-\bar{\gamma_{12}} & \bar{\gamma_{22}} & -\bar{\gamma_{32}}\\
\end{array} } \right).
\]
Substituting $\tilde{Q}\left( z \right)-S$, $\left( \tilde{A}-z 
\right)^{-1}$ and $\tilde{\Gamma} $ into (\ref{eq104}), we obtain: 
\[
\left( {\begin{array}{*{20}c}
\frac{1}{z\left( z+1 \right)} & \frac{-2}{z+1}\\
\frac{-2}{z+1} & 0\\
\end{array} } \right)=
\]
\[
=\left( {\begin{array}{*{20}c}
-\bar{\gamma_{11}} & \bar{\gamma_{21}} & -\bar{\gamma_{21}}\\
-\bar{\gamma_{12}} & \bar{\gamma_{22}} & -\bar{\gamma_{32}}\\
\end{array} } \right)\left( {\begin{array}{*{20}c}
-\frac{1}{z} & 0 & 0\\
0 & \frac{-1}{z+1} & 0\\
0 & 0 & \frac{-1}{z+1}\\
\end{array} } \right)\left( {\begin{array}{*{20}c}
\gamma_{11} & \gamma_{12}\\
\gamma_{21} & \gamma_{22}\\
\gamma_{31} & \gamma_{32}\\
\end{array} } \right).
\]
Solving this nonlinear system for the unknowns $\gamma_{ij},\, 
i=1,2,3;j=1,2$, yields:
\[
\tilde{Q}\left( z \right)=\left( {\begin{array}{*{20}c}
0 & 1\\
1 & 0\\
\end{array} } \right)+\left( {\begin{array}{*{20}c}
-1 & 1 & 0\\
0 & \-2 & -2\\
\end{array} } \right)\left( {\begin{array}{*{20}c}
-\frac{1}{z} & 0 & 0\\
0 & \frac{-1}{z+1} & 0\\
0 & 0 & \frac{-1}{z+1}\\
\end{array} } \right)\left( {\begin{array}{*{20}c}
1 & 0\\
1 & 2\\
0 & 2\\
\end{array} } \right).
\]
Finally, by multiplying this expression for $\tilde{Q}\left( z \right)$ by $(z+1)$, we obtain the representation (\ref{eq108}) of $Q$.\hfill $\square$ 

\[
Q\left( z \right)=\left( {\begin{array}{*{20}c}\frac{1}{z^{2}} & 1-\frac{1}{z}\\1-\frac{1}{z} & 0\\
\end{array} } \right)
\]

\end{document}